\documentclass[preprint,authoryear,11pt]{elsarticle}

\usepackage{amsthm,amsmath,amssymb,mathrsfs,ulem,graphicx}
\newtheorem{definition}{Definition}
\newtheorem{lemma}{Lemma}
\newtheorem{proposition}{Proposition}
\newtheorem{theorem}{Theorem}
\newtheorem{corollary}{Corollary}

\newcommand{\dif}{\mathrm{d}}
\newcommand{\0}{\mathbf{0}}
\newcommand{\bfalpha}{\boldsymbol{\alpha}}
\newcommand{\bftheta}{\boldsymbol{\theta}}

\newcommand{\bfbeta}{\boldsymbol{\beta}}

\newcommand{\bftau}{\boldsymbol{\tau}}
\newcommand{\bfeta}{\boldsymbol{\eta}}
\newcommand{\bfSigma}{\boldsymbol{\Sigma}}
\newcommand{\bfrmA}{\mathbf{A}}
\newcommand{\bfrmB}{\mathbf{B}}
\newcommand{\bfc}{\boldsymbol{c}}
\newcommand{\bfh}{\boldsymbol{h}}
\newcommand{\bfs}{\boldsymbol{s}}
\newcommand{\bfx}{\boldsymbol{x}}
\newcommand{\bfy}{\boldsymbol{y}}
\newcommand{\bfX}{\boldsymbol{X}}
\newcommand{\Euler}{\mathrm{e}}
\newcommand{\Reals}{\mathbb{R}}

\newcommand{\PSD}{\mathbb{S}_{+}}
\newcommand{\calT}{\mathcal{T}}
\newcommand{\calD}{\mathcal{D}}
\newcommand{\T}{{}^{\mathsf{T}}}
\newcommand{\BoltzmanK}{k_{\mathrm{B}}}
\newcommand{\E}{\mathcal{E}}
\newcommand{\Var}{\mathcal{V}}
\newcommand{\Cov}{\mathcal{C}}
\newcommand{\Normal}{\mathscr{N}}
\newcommand{\FIM}{\mathscr{I}}
\newcommand{\TheRest}{\mathscr{O}}
\newcommand{\trace}{\mathrm{tr}}
\newcommand{\stdefficiency}{\mathrm{eff}}
\newcommand{\efficiency}{\mathrm{ueff}}
\newcommand{\ord}{\mathrm{o}}
\newcommand{\romanI}{\mathrm{I}}
\newcommand{\romanII}{\mathrm{II}}
\newcommand{\OUPClass}{\mathcal{OU}}
\newcommand{\LoewnerDominates}{\succeq_{\mathrm{L}}}
\newcommand{\LoewnerDominatedBy}{\preceq_{\mathrm{L}}}

\biboptions{longnamesfirst,comma}

\journal{Journal of Statistical Planning and Inference}

\begin{document}

\begin{frontmatter}



\title{Ultimate efficiency of designs for processes of Ornstein-Uhlenbeck type}


\author{Vladim\'{i}r Lacko}
\ead{lackovladimir@gmail.com}
\address{Department of Applied Mathematics and Statistics, Faculty of Mathematics, Physics and Informatics, Comenius University, 842 48 Bratislava, Slovak Republic}

\begin{abstract}
For a process governed by a linear It\=o stochastic differential equation of the form $\dif X(t) = [a(t)+b(t)X(t)]\dif t + \sigma(t)\dif W(t)$ we prove an existence of optimal sampling designs with strictly increasing sampling times. We derive an asymptotic Fisher information matrix, which we take as a reference in assessing a quality of finite-point sampling designs. The results are extended to a broader class of It\=o stochastic differential equations satisfying a certain condition. We give an example based on the Gompertz growth law refuting a generally accepted opinion that small-sample designs lead to a very high level of efficiency.
\end{abstract}


\end{frontmatter}

\section{Introduction}

Suppose we can observe a univariate continuous-time process $\{X(t)\}_{t\geq0}$ governed by a linear It\=o stochastic differential equation of the form
\begin{eqnarray}
\dif X(t) &=& [a_{\bftheta_1,\bftheta_2,\theta_3}(t)+b_{\bftheta_2,\theta_3}(t)X(t)]\dif t + \sigma_{\theta_3}(t)\dif W(t)\label{eqn:MODEL}\\
&=& f_{\bftheta}(t,X(t))\dif t + \sigma_{\theta_3}(t)\dif W(t),\nonumber\\
X(0) &=& X_0\in\Reals\;\text{is fixed},\nonumber
\end{eqnarray}
at $n$ strictly increasing design times $\bftau=(t_1,\ldots,t_n)$ from the set of feasible $n$-point (sampling) designs
\begin{equation*}
\calT_{n,\calD}=\{\bftau\in\Reals^n:T_*\leq t_1<\cdots<t_n\leq T^*\}
\end{equation*}
on the experimental domain $\calD=[T_*,T^*]$, where $0<T_*<T^*$. Here, the function $f(t,x)$ is the drift, $\sigma(t)$ is the volatility and $\{W(t)\}_{t\geq0}$ is a one-dimensional Wiener process.

The observed quantities $X(t_1),\ldots,X(t_n)$ depend on unknown vector parameters $\bftheta_1\in\Reals^{m_1}$, $\bftheta_2\in\Reals^{m_2}$ and a scalar parameter $\theta_3\in\Reals$, which appear in the coefficients of equation \eqref{eqn:MODEL}; where needed we emphasize this dependence by the corresponding subscripts. The fixed initial value $X_0\in\Reals$, which absents in the governing equation, might be in some instances regarded as an unknown parameter, too. One may wonder about the choice of the parametrisation. In Section \ref{sec:NONLINEAR_MODEL_FOR_OBSERVATIONS}, where we formulate a nonlinear regression model for observations, we show that $X_0$ and $\bftheta_1$ are present only in the mean value, while $\bftheta_2$ is a parameter common for the mean value and the covariance matrix of observations, and $\theta_3$ is a variance parameter but with asymptotic properties different to those for $\bftheta_2$; we give a detailed discussion in Section \ref{sec:ULTIMATE_EFFICIENCY_OF_DESIGNS}. Such discrimination enables us to understand what roles particular parameter groups play in the model. For the sake of simplicity we define $\bftheta=(X_0,\bftheta_1\T,\bftheta_2\T,\theta_3)\T$ and assume that $n\geq\dim(\bftheta)=m$. We further require the functions $a(t)$, $b(t)$ and $\sigma^2(t)$ to be differentiable with respect to $\bftheta$ and integrable with respect to $t$, and $\sigma(t)$ to be positive almost everywhere with respect to Lebesgue measure on the real axis.

Searching for feasible optimal designs, whose existence we prove in Section \ref{sec:FISHER_INFORMATION_MATRIX_AND_EXISTENCE_OF_FEASIBLE_LOCALLY_OPTIMAL_DESIGNS}, might be computationally challenging. Instead we therefore try to answer a different but related question: how can we evaluate a quality of sampling designs for maximum likelihood estimation of the parameter subvector $(X_0,\bftheta_1\T,\bftheta_2\T)\T$? One option is given in Section \ref{sec:ULTIMATE_EFFICIENCY_OF_DESIGNS}, where as a benchmark we take the asymptotic Fisher information matrix, which results from the observation of the full trajectory.

Governing equation \eqref{eqn:MODEL} is motivated by a variety of problems studied in literature. The first group of problems consists of modifications of Ornstein-Uhlenbeck process. \cite{UHLENBECK_ORNSTEIN_1930} proposed a simple model of particle velocity $X(t)$ that can be rewritten to a stochastic differential equation
\begin{equation}\label{eqn:ORNSTEIN_UHLENBECK_PROCESS}
\dif X(t)=\left(\frac{\theta_1}{\theta_2}-\frac{\theta_3}{\theta_2} X(t)\right)\dif t + \left(\frac{2\BoltzmanK T}{\theta_3}\right)^{1/2}\dif W(t),\;X(0)=X_0\;\text{fixed},
\end{equation}
where $\theta_1$ is the level of external force, $\theta_2$ is the mass, $\theta_3$ is the friction coefficient, $T$ is the temperature of the system and $\BoltzmanK$ is the Boltzmann constant. In applications we can find a more abstract arrangement of stochastic differential equation \eqref{eqn:ORNSTEIN_UHLENBECK_PROCESS} given by
\begin{equation}\label{eqn:SIMPLE_ORNSTEIN_UHLENBECK_PROCESS}
\dif X(t) = \theta_2[\theta_1 - X(t)]\dif t + \theta_3\dif W(t),\;X(0)=X_0\;\text{fixed},
\end{equation}
where the asymptotic mean $\theta_1$, the mean-reversion speed $\theta_2$ and the diffusion coefficient $\theta_3$ are the unknown parameters. For instance, \cite{RICCARDI_SACERDOTE_1979} used the process governed by equation \eqref{eqn:SIMPLE_ORNSTEIN_UHLENBECK_PROCESS} to model the voltage difference between the membrane and resting potentials at the trigger zone of the neuron. Note that compared to equation \eqref{eqn:ORNSTEIN_UHLENBECK_PROCESS} the drift and diffusion part of the simpler variant \eqref{eqn:SIMPLE_ORNSTEIN_UHLENBECK_PROCESS} do not have any common parameters.

Applications of Ornstein-Uhlenbeck process lead to investigations of its optimal sampling. For a process driven by stochastic differential equation \eqref{eqn:SIMPLE_ORNSTEIN_UHLENBECK_PROCESS} with $\theta_1$ and $X_0$ regarded as unknown parameters, \cite{HARMAN_STULAJTER_2009} proved optimality of equidistant sampling designs for estimation and prediction. Studies carried out by \cite{KISELAK_STEHLIK_2008} and later \cite{ZAGORAIOU_ANTOGNINI_2009} show that equidistant designs are optimal also for a stationary variant of equation \eqref{eqn:SIMPLE_ORNSTEIN_UHLENBECK_PROCESS}, that is, for $t\to\infty$. Further analysis of model \eqref{eqn:SIMPLE_ORNSTEIN_UHLENBECK_PROCESS} with $\theta_1$ and $X_0$ in the position of parameters was provided by \cite{LACKO_2012}, who demonstrated that under the assumption of a time-dependent volatility the optimal sampling times are more concentrated in areas with lower levels of the volatility function. 

The second group of problems covered by model \eqref{eqn:MODEL} is a family of Brownian motions of the form $\{A_{\bftheta_1}(t)+\theta_3 W_t\}_{t\geq0}$ with unknown parameters $\bftheta_1$ and $\theta_3$, which coincides with the stochastic differential equation
\begin{equation*}
\dif X(t) = a_{\bftheta_1}(t)\dif t + \theta_3\dif W(t),\;X(0)=0,
\end{equation*}
where $a_{\bftheta_1}(t)$ is the derivative of $A_{\bftheta_1}(t)$ with respect to $t$.

Brownian motions were to a certain extend analysed by \cite{SACKS_YLVISAKER_1966,SACKS_YLVISAKER_1968}, who proposed an approach for construction of the so-called asymptotically optimal designs that work well for a large number of observations. On the contrary, for a Brownian motion with $A_{\bftheta_1}(t)=(1,t,t^2)\bftheta_1$, \cite{HARMAN_STULAJTER_2011} showed that exact equidistant sampling designs are optimal for parameter estimation as well as prediction.


Throughout this paper we use the following {\it notation:} For any random or non-random function $h(t)$ and design $\bftau\in\calT_{n,\calD}$ we denote $\bfh(\bftau)=(h(t_1),\ldots,h(t_n))\T$. Similarly, for a symmetric function $\Sigma(t_i,t_j)=\Sigma(t_j,t_i)$ we have $\{\bfSigma(\bftau)\}_{ij} = \Sigma(t_i,t_j)$, $1\leq i,j\leq n$. The norm of a design $\bftau\in\calT_{n,\calD}$ is defined as $\|\bftau\|=\max_{2\leq i\leq n}(t_i-t_{i-1})$. Symbols $\E$, $\Var$ and $\Cov$ represent expected value (or vector), variance (variance-covariance matrix) and covariance, respectively. By $\PSD^m$ we denote the set of all $m\times m$ symmetric non-negative definite matrices. We say that matrix $\bfrmA$ ``Loewner dominates'' matrix $\bfrmB$, denoted by $\bfrmA\LoewnerDominates\bfrmB$, if $\bfrmA-\bfrmB\in\PSD^m$.

\section{Nonlinear model for observations}\label{sec:NONLINEAR_MODEL_FOR_OBSERVATIONS}

To evaluate the Fisher information matrix for a given design $\bftau\in\calT_{n,\calD}$ we need to understand mutual relations between individual components of $\bfX(\bftau)=(X(t_1),\ldots,X(t_n))\T$.

It is a well-known fact that by applying It\=o's lemma (see Theorem~6 by \cite{ITO_1951}) to transformation $\Euler^{-B(t)}X(t)$, where $B(t)$ is an arbitrary antiderivative of $b(t)$, we obtain $X(t)=\E[X(t)] + \varepsilon(t)$ for all $t\geq0$, where
\begin{equation}\label{eqn:MEAN_VALUE}
\E[X(t)] = \Euler^{B(t)-B(0)}X_0 + \int_0^t\Euler^{B(t)-B(\nu)}a(\nu)\dif\nu
\end{equation}
is the expectation of the process $\{X(t)\}_{t\geq0}$ at time $t$ and
\begin{equation*}
\varepsilon(t) = \int_0^t\Euler^{B(t)-B(\nu)}\sigma(\nu)\dif W(\nu)
\end{equation*}
is a zero-mean Gaussian random variable; see, e.g., \cite{GARDINER_1985}.

\begin{lemma}\label{lem:COVARIANCE}
Let $\{X(t)\}_{t\geq0}$ be a process governed by stochastic differential equation \eqref{eqn:MODEL}. Then
\begin{equation*}
\Var[X(t)] = \int_0^t \Euler^{2[B(t)-B(\nu)]}\sigma^2(\nu)\dif\nu
\end{equation*}
and
\begin{equation*}
\forall_{0\leq t_1\leq t_2}\;\Cov[X(t_1),X(t_2)] = \Euler^{B(t_2)-B(t_1)}\Var[X(t_1)].
\end{equation*}
\end{lemma}

\begin{proof}
The expression for $\Var[X(t)]$ comes from It\=o's isometry; see \cite{OKSENDAL_2000}. To show the second part of the lemma we make out an ordinary differential equation for the covariance function. Since $\{X(t)\}_{t>0}$ is a Gaussian process, usual regularity conditions for interchange of differentiation and integration order are satisfied. Taking the expectation of the governing equation \eqref{eqn:MODEL} yields $\frac{\dif}{\dif s}\E[X(t+s)]=\E[f(t+s,X(t+s))]$, where we recall that $f(t,x)=a(t)+b(t)x$. Consequently, for all $t,s\geq0$
\begin{eqnarray*}
\frac{\dif}{\dif s}\E[X(t)X(t+s)] &=& \E\left[X(t)\frac{\dif}{\dif s}\E[X(t+s)]\right] = \E[X(t)f(t+s,X(t+s))],\\
\E[X(t)]\frac{\dif}{\dif s}\E[X(t+s)] &=& \E[X(t)]\E[f(t+s,X(t+s))].
\end{eqnarray*}
Basic rules for covariance give an ordinary differential equation
\begin{equation*}
\frac{\dif}{\dif s}\Cov[X(t),X(t+s)] = \Cov[X(t),f(t+s,X(t+s))] = b(t+s)\Cov[X(t),X(t+s)],
\end{equation*}
with the initial value $\Cov[X(t),X(t+s)]=\Var[X(t)]$ at $s=0$. A use of standard methods for solving ordinary differential equations and subsequent setting $t_1=t$ and $t_2=t+s$ entail the statement of the lemma.
\end{proof}

\begin{proposition}\label{prop:COVARIANCE_STRUCTURE}
Let $\{X(t)\}_{t\geq0}$ be a process governed by stochastic differential equation \eqref{eqn:MODEL} and let $\bftau\in\calT_{n,\calD}$ be a given design. Then the $ij$th component, $0\leq i\leq j\leq n$, of variance-covariance matrix $\Var[\bfX(\bftau)]=\bfSigma(\bftau)$ is
\begin{eqnarray*}
\{\bfSigma(\bftau)\}_{ij} &=& u(t_i)v(t_j),\;\text{where}\\
u(t) &=& \Euler^{B(t)}\int_0^t\Euler^{-2B(\nu)}\sigma^2(\nu)\dif\nu,\\
v(t) &=& \Euler^{B(t)}.
\end{eqnarray*}
Moreover, for any $\bftau\in\calT_{n,\calD}$ the matrix $\bfSigma(\bftau)$ is positive definite.
\end{proposition}

\begin{proof}
The form of the covariance kernel comes from Lemma~\ref{lem:COVARIANCE}. To prove the positive definiteness of $\bfSigma(\bftau)$ we use the induction. Notice that for any continuous function $g(t)$ positive almost everywhere with respect to the Lebesgue measure on the real axis we have
\begin{equation}\label{eqn:VARIANCE_INEQUALITY}
D(t_2) = \Euler^{2B(t_2)}\int_0^{t_2}\Euler^{-2B(\nu)}g(\nu)\dif\nu>\Euler^{2[B(t_2)-B(t_1)]}D(t_1),
\end{equation}
where $D(t)$ coincides with $\Var[X(t)]$ for $g(t)=\sigma^2(t)$. We can use Lemma~\ref{lem:COVARIANCE} and inequality \eqref{eqn:VARIANCE_INEQUALITY} to show that for any feasible $2$-point design the covariance matrix is positive definite. Let us assume that for $\bftau_n\in\calT_{n,\calD}$ the matrix $\bfSigma(\bftau_n)$ is positive definite and, without loss of generality, let $\bftau_{n+1}=(\bftau_n\T,t_{n+1})\T\in\calT_{n+1,\calD}$. The matrix $\bfSigma(\bftau_{n+1})$ is row-equivalent to
\begin{equation*}
\begin{pmatrix}
\bfSigma(\bftau_n) & \bfs\\
\0_m\T & \Var[X(t_{n+1})] - \bfs\T\bfSigma(\bftau_n)\bfs
\end{pmatrix},
\end{equation*}
where $\bfs=(\Euler^{B(t_{n+1})-B(t_1)}\Var[X(t_1)],\ldots,\Euler^{B(t_{n+1})-B(t_n)}\Var[X(t_n)])\T$. By using the technique of \cite{HARMAN_STULAJTER_2009} we can show that
\begin{equation*}
\Var[X(t_{n+1})] - \bfs\T\bfSigma(\bftau_n)\bfs = v^2(t_{n+1})\left(\frac{u(t_{n+1})}{v(t_{n+1})}-\frac{u(t_n)}{v(t_n)}\right) = \Euler^{2B(t_{n+1})}\int_{t_n}^{t_{n+1}}\Euler^{-2B(\nu)}g(\nu)\dif\nu > 0,
\end{equation*}
where we adopted the notation from inequality \eqref{eqn:VARIANCE_INEQUALITY}.
\end{proof}


The expectation of the underlying model \eqref{eqn:MEAN_VALUE} and Proposition \ref{prop:COVARIANCE_STRUCTURE} enable us to formulate the design problem for stochastic differential equation \eqref{eqn:MODEL} in terms of nonlinear regression
\begin{equation}\label{eqn:REGRESSION_MODEL}
\forall_{\bftau\in\calT_{n,\calD}}\;\bfX(\bftau)\sim\Normal\Big(\E_{X_0,\bftheta_1,\bftheta_2,\theta_3}[\bfX(\bftau)],\bfSigma_{\bftheta_2,\theta_3}(\bftau)\Big),
\end{equation}
where $\E[\bfX(\bftau)]$ depends on all unknown parameters, while $\bfSigma(\bftau)$ is influenced only by $\bftheta_2$ and $\theta_3$.

\section{Fisher information matrix and existence of feasible locally optimal designs}\label{sec:FISHER_INFORMATION_MATRIX_AND_EXISTENCE_OF_FEASIBLE_LOCALLY_OPTIMAL_DESIGNS}

For maximum likelihood estimation the Fisher information matrix is a standard reference for a quality of estimated values of unknown parameters, and experimenter naturally prefers designs with ``large'' Fisher information matrix.

Discrimination between two designs in a single-parameter model is naive since the Fisher information matrix is reduced to a scalar value. On the other hand, for multi-parameter setup, where matrices need to be compared, we measure the ``size'' of information matrices by information functions. An information function is any function $\Phi:\PSD^m\mapsto[0,\infty)$ which is non-constant, concave, upper semi-continuous, positive homogeneous and Loewner isotonic (i.e., if  $\FIM_1\LoewnerDominates\FIM_2$ then $\Phi[\FIM_1]\geq\Phi[\FIM_2]$). 
The most popular information functions, such as those corresponding to D-, E-, A- and $\bfc$-optimality, conform to reasonable geometrical or statistical criteria. We refer the reader to \cite{PAZMAN_1986} and \cite{PUKELSHEIM_1993} for more details on information functions and optimality criteria.

Under regression model \eqref{eqn:REGRESSION_MODEL} particular blocks of the Fisher information matrix take the form
\begin{equation*}
\{\FIM(\bftau,\bftheta^*)\}_{\bfalpha\bfbeta} = \left(\frac{\partial\E[\bfX(\bftau)]}{\partial \bfalpha}\bfSigma^{-1}(\bftau)\frac{\partial\E[\bfX(\bftau)]}{\partial \bfbeta\T}\right)\Bigg|_{\bftheta^*} + \frac{1}{2}\trace\left\{\bfSigma^{-1}(\bftau)\frac{\partial \bfSigma(\bftau)}{\partial \bfalpha}\bfSigma^{-1}(\bftau)\frac{\partial \bfSigma(\bftau)}{\partial \bfbeta\T}\right\}\Bigg|_{\bftheta^*},
\end{equation*}
where $\bfalpha,\bfbeta=X_0,\bftheta_1,\bftheta_2,\theta_3$, and $\bftheta^*$ is a guess at the true value of unknown vector parameter, see \cite{MARDIA_MARSHALL_1984}.

\begin{theorem}\label{thm:OPTIMAL_DESIGN_FOR_INITIAL_VALUE}
If the initial value $X_0$ of stochastic differential equation \eqref{eqn:MODEL} is the only unknown parameter, then it is optimal to take $t_1=T_*$ regardless of the number of design points. The variance of the corresponding maximum likelihood estimate is
\begin{equation*}
\Var[\hat{X}_0] = \int_0^{T_*}\Euler^{2[B(0)-B(\nu)]}\sigma^2(\nu)\dif\nu.
\end{equation*}
\end{theorem}

\begin{proof}
The key for the proof is the result of \cite{HARMAN_STULAJTER_2009}: for a positive definite matrix $\bfrmA$ such that $\{\bfrmA\}_{ij}=u_i v_j$, $1\leq i\leq j\leq n$, we have
\begin{equation}\label{eqn:HARMAN_STULAJTER_FORMULA}
\forall_{\bfx,\bfy\in\Reals^n}\;\bfx\T\bfrmA^{-1}\bfy = \frac{x_1y_1}{u_1v_1} + \sum_{i=2}^n\frac{\left(\frac{x_i}{v_i}-\frac{x_{i-1}}{v_{i-1}}\right)\left(\frac{y_i}{v_i}-\frac{y_{i-1}}{v_{i-1}}\right)}{\frac{u_i}{v_i}-\frac{u_{i-1}}{v_{i-1}}}.
\end{equation}
Since for any $X_0^*$
\begin{equation*}
\left.\frac{\partial\E[X(t)]}{\partial X_0}\right|_{X_0^*} = v(t)\Euler^{-B(0)} \;\text{and}\;\left.\frac{\partial\Var[X(t)]}{\partial X_0}\right|_{X_0^*} = 0,
\end{equation*}
from relation \eqref{eqn:HARMAN_STULAJTER_FORMULA} and Proposition \ref{prop:COVARIANCE_STRUCTURE} we obtain that for any design $\bftau\in\calT_{n,\calD}$ the information about $X_0$ is
\begin{equation*}
\FIM(\bftau,X_0^*) = \Var^{-1}[\hat{X}_0] = \left(\int_0^{t_1}\Euler^{2[B(0)-B(\nu)]}\sigma^2(\nu)\dif\nu\right)^{-1},
\end{equation*}
that is, the information is determined by the first sampling time $t_1$. The function $\Euler^{2[B(0)-B(t)]}\sigma^2(t)$ is positive for any $t$, therefore $\FIM(\bftau,X_0^*)$ is maximal for $t_1=T_*$.
\end{proof}

Although the statement of Theorem \ref{thm:OPTIMAL_DESIGN_FOR_INITIAL_VALUE} might give an impression to hold true for stochastic processes in general, the opposite is true. For instance, consider a process $\{X(t)\}_{t\geq0}$ governed by equation
\begin{equation}\label{eqn:COUNTEREXAMPLE}
\dif X(t) = \frac{\theta_2}{2}(X_0 - X(t))\dif t + \Euler^{-\theta_3 t/2}\dif W(t),\; X(0)=X_0,
\end{equation}
for which we can perform only one observation. The values $\theta_2>0$ and $\theta_3>0$ are known, while $X_0$ is an unknown parameter. Evidently, model \eqref{eqn:COUNTEREXAMPLE} violates the assumption on absence of the initial value in the governing equation. For this model the mean value is equal to $X_0$ for all $t\geq0$ and the Fisher information for model \eqref{eqn:COUNTEREXAMPLE} obtained from the observation performed at the time $t$, which equals to the reciprocal of the variance, attains the value
\begin{equation*}
\FIM(t,\theta_2,\theta_3) = \left\{\begin{array}{rl}
\dfrac{\theta_2-\theta_3}{\Euler^{-\theta_3 t}-\Euler^{-\theta_2 t}},&\theta_2\neq\theta_3,\\
\dfrac{1}{\Euler^{-\theta_2 t}t},&\theta_2=\theta_3.
\end{array}\right.
\end{equation*}
Since $\FIM(t,\theta_2,\theta_3)$ a convex function, for $\theta_2\neq\theta_3$ we can find that the observation time minimising the Fisher information is given by $t_{\min \FIM}=(\ln \theta_2 - \ln \theta_3)/(\theta_2 - \theta_3)$ and $t_{\min \FIM}=1/\theta_2$ as $\theta_3\to\theta_2$. Consequently, if the bounds of the experimental domain satisfy an inequality $T_*<t_{\min \FIM}<T^*$ with $\FIM(T_*,\theta_2,\theta_3)<\FIM(T^*,\theta_2,\theta_3)$ or the bounds satisfy $t_{\min \FIM}\leq T_*<T^*$ then it is optimal to observe the process as late as posible.

Under a more general setup, where the parametrisation of the model is not as simple as in Theorem~\ref{thm:OPTIMAL_DESIGN_FOR_INITIAL_VALUE}, to find an optimal design we usually employ a battery of optimization procedures. The set of feasible designs $\calT_{n,\calD}$, $n\geq m=\dim(\bftheta)$, is not closed and optimization algorithms might tend to designs, which do not meet the requirements of the experimental setup. Therefore, it is important to ensure an existence of a feasible optimal $n$-point design, that is, belonging to $\calT_{n,\calD}$.

\begin{lemma}\label{lem:INFORMATION_FOR_MARKOV_PROCESS}
Let $\{X(t)\}_{t\geq0}$, $X(0)$ fixed, be a $\bftheta$-parametrised continuous-time Markov process with transition density kernel satisfying usual regularity conditions; see (5.12) in \cite{LEHMANN_CASELLA_1998}. Then for any $\bftau\in\calT_{n,\calD}$ the Fisher information matrix for $\bfX(\bftau)$ is
\begin{equation*}
\FIM(\bftau,\bftheta^*) = \FIM_{X(t_1)\mid X(0)}(\bftheta^*) + \sum_{i=2}^n\E_{X(t_{i-1})}\left[\FIM_{X(t_i)\mid X(t_{i-1})}(\bftheta^*)\right],
\end{equation*}
where $\FIM_{X(t_i)\mid X(t_{i-1})}(\bftheta^*)$ is the Fisher information matrix for $X(t_i)$ conditioned on the value of $X(t_{i-1})$ and $\E_{X(t_{i-1})}[\cdot]$ is the expectation with respect to $X(t_{i-1})$.
\end{lemma}

\begin{proof}
We denote $t_0=0$ and consider the derivatives to be evaluated at $\bftheta^*$. Let $p(x,s\mid y,t)=\frac{\dif}{\dif x}\Pr[X(t+s)<x\mid X(t)=y]$ be the transition density kernel of the process $\{X(t)\}_{t\geq0}$, where $p(x,0\mid y,t)=\delta(x-y)$ and $\delta(\cdot)$ is the Dirac delta function. Then
\begin{eqnarray*}
\FIM(\bftau,\bftheta^*) &=& \E_{\bfX(\bftau)}\left[\frac{\partial^2}{\partial\bftheta\partial\bftheta\T}\ln\left(\prod_{i=1}^np(x_i,t_i-t_{i-1}\mid x_{i-1},t_{i-1})\right)\right]\\
&=& \int_{x_1\in\Reals}\cdots\int_{x_n\in\Reals}\sum_{i=1}^n\left(\frac{\partial^2}{\partial\bftheta\partial\bftheta\T}\ln p(x_i,t_i-t_{i-1}\mid x_{i-1},t_{i-1})\right)\\
&&\times\prod_{k=1}^np(x_k,t_k-t_{k-1}\mid x_{k-1},t_{k-1})\dif x_1\ldots\dif x_n\\
&=& \sum_{i=1}^n\int_{x_1\in\Reals}\cdots\int_{x_n\in\Reals}\frac{\partial^2}{\partial\bftheta\partial\bftheta\T}\ln p(x_i,t_i-t_{i-1}\mid x_{i-1},t_{i-1})\\
&&\times\prod_{k=1}^np(x_k,t_k-t_{k-1}\mid x_{k-1},t_{k-1})\dif x_1\ldots\dif x_n\\
&=& \sum_{i=1}^n\int_{x_{i-1}\in\Reals}\int_{x_i\in\Reals}\left(\frac{\partial^2}{\partial\bftheta\partial\bftheta\T}\ln p(x_i,t_i-t_{i-1}\mid x_{i-1},t_{i-1})\right)\\
&&\times p(x_i,t_i-t_{i-1}\mid x_{i-1},t_{i-1})p(x_{i-1},t_{i-1}\mid X(0),0)\dif x_{i-1}\dif x_i\\
&=& \sum_{i=1}^n\int_{x_{i-1}\in\Reals}\FIM_{X(t_i)\mid X(t_{i-1})}(\bftheta^*)p(x_{i-1},t_{i-1}\mid X(0),0)\dif x_{i-1}.
\end{eqnarray*}
\end{proof}

\begin{theorem}\label{thm:EXISTENCE_OF_NONDEGENERATE_OPTIMAL_DESIGN}
Let $\{X(t)\}_{t\geq0}$ be a process governed by stochastic differential equation \eqref{eqn:MODEL} and $\overline{\calT_{n,\calD}}$ be the closure of $\calT_{n,\calD}$. If $\bftau_0\in\overline{\calT_{n,\calD}}$ then there exists a feasible design $\bftau\in\calT_{n,\calD}$ such that $\FIM(\bftau,\bftheta^*)\LoewnerDominates\FIM(\bftau_0,\bftheta^*)$.

\end{theorem}

\begin{proof}
If $\bftau_0\in\calT_{n,\calD}$ then $\bftau=\bftau_0$. Assume that $\bftau_0\in\overline{\calT_{n,\calD}}\setminus\calT_{n,\calD}$, i.e., there exists $i\in\{2,\ldots,n\}$ such that $t_{i-1}=t_i$. Equation \eqref{eqn:MEAN_VALUE} and Proposition \ref{prop:COVARIANCE_STRUCTURE} yield $X(t)\mid X(t) = X(t)$ almost surely, which gives no information about unknown parameter, so $\E_{X(t_i)}[\FIM_{X(t_i)\mid X(t_i)}(\bftheta^*)] = \0_{m\times m}$. As a consequence of Lemma~\ref{lem:INFORMATION_FOR_MARKOV_PROCESS}, by leaving $t_i$ from the experimental design the amount of information does not change. We repeat this procedure until we obtain a design $\bftau_1\in\calT_{n_1,\calD}$, $n_1<n$, for which we have $\FIM(\bftau_1,\bftheta^*)=\FIM(\bftau_0,\bftheta^*)$. For $\bftau$ we take any design from $\calT_{n,\calD}$ with arbitrary $n_1$ components being given by $\bftau_1$. Analogously to Lemma~\ref{lem:INFORMATION_FOR_MARKOV_PROCESS} we can show that $\FIM(\bftau,\bftheta^*) = \FIM_{\bfX(\bftau_1)\mid X(0)}(\bftheta^*) + \E_{\bfX(\bftau_1)}[\FIM_{\bfX(\bftau_2)\mid \bfX(\bftau_1)}(\bftheta^*)] = \FIM(\bftau_1,\bftheta^*) + \E_{\bfX(\bftau_1)}[\FIM_{\bfX(\bftau_2)\mid \bfX(\bftau_1)}(\bftheta^*)]$, which implies
\begin{equation*}
\FIM(\bftau,\bftheta^*) - \FIM(\bftau_0,\bftheta^*) = \E_{\bfX(\bftau_1)}[\FIM_{\bfX(\bftau_2)\mid \bfX(\bftau_1)}(\bftheta^*)]\LoewnerDominates\0_{m\times m}.
\end{equation*}
\end{proof}

Obviously, statement of Theorem~\ref{thm:EXISTENCE_OF_NONDEGENERATE_OPTIMAL_DESIGN} holds true also in more general models with correlated observations: if $\E_{X(t)}[\FIM_{X(t)\mid X(t)}(\bftheta^*)] = \0_{m\times m}$ for all $t\in\calD$ then there exists a feasible optimal design. Additionally, feasible designs in accordance with Theorem~\ref{thm:EXISTENCE_OF_NONDEGENERATE_OPTIMAL_DESIGN} dominate ``degenerated'' designs with respect to Loewner ordering, thus, feasible designs are universally as good as ``degenerated'' or better.

Note that the conditioned Fisher information $\FIM_{\bfX(t+s)\mid \bfX(t)}(\bftheta^*)$ is not continuous for $s$ approaching $0$. Under model \eqref{eqn:MODEL} it follows from relation \eqref{eqn:CONDITIONED_VARIANCE_EXPANSION} given later that $\lim_{s\to\infty}\{\FIM_{\bfX(t+s)\mid \bfX(t)}(\bftheta^*)\}_{\theta_3\theta_3}\neq \{\FIM_{\bfX(t)\mid \bfX(t)}(\bftheta^*)\}_{\theta_3\theta_3}$. In the last two paragraphs in the next section we discuss this fact in connection with the consistency of maximum likelihood estimator of $\theta_3$.

\section{Ultimate efficiency of designs}\label{sec:ULTIMATE_EFFICIENCY_OF_DESIGNS}

In many situations we do not need to use optimal designs but we accept also designs that are in a certain way relatively efficient. A standard approach to measuring relative efficiency is comparison of a given design $\bftau\in\calT_{n,\calD}$ with an optimal design by considering
\begin{equation}\label{eqn:EFFICIENCY}
\stdefficiency(\bftau\mid\Phi,\bftheta^*) = \frac{\Phi[\FIM(\bftau,\bftheta^*)]}{\sup_{\bfeta\in\calT_{n,\calD}}\Phi[\FIM(\bfeta,\bftheta^*)]},
\end{equation}
see \cite{PUKELSHEIM_1993}. As an alternative to evaluating the optimal information function in the denominator of ratio \eqref{eqn:EFFICIENCY}, instead of $\sup_{\bfeta\in\calT_{n,\calD}}\Phi[\FIM(\bfeta,\bftheta^*)]$ \cite{PAZMAN_2007} suggested to take a value $\Phi[\tilde{\FIM}]$, where $\tilde{\FIM}$ is any suitable reference matrix. Particularly, under standard regression model with correlated observations, where we are interested in the parameters of the mean value which do not to appear in the covariance structure, \cite{PAZMAN_2007} further noted that a suitable choice for the reference matrix $\tilde{\FIM}$ is the asymptotic Fisher information matrix obtained by observation of the full trajectory. Such asymptotic Fisher information matrix has all eigenvalues finite, that is, the maximum likelihood estimator is not consistent. Analogously, the value of the corresponding information function attains a finite value as well. The idea to measure the relative efficiency of a design with respect to information obtained by observation of the full trajectory, which we call ``ultimate efficiency'' as \cite{HARMAN_2011} suggested, was later adopted in the works of \cite{HARMAN_STULAJTER_2011} and \cite{LACKO_2012}, and plays a crucial role also in this paper. 

Compared with the standard setup the parametrisation of the model \eqref{eqn:MODEL} has, however, a different structure. In a more general situation, if we observe the full trajectory, i.e., we perform measurements at design points $\bftau^{(n)}=(T_*,t_2^{(n)},\ldots,t_{n-1}^{(n)},T^*)\T\in\calT_{n,\calD}$ on a bounded domain $\calD$ and $\|\bftau^{(n)}\|\to0$ as $n\to\infty$, then the Fisher information matrix $\FIM(\bftau^{(n)},\bftheta^*)$ converges to a matrix $\FIM_{\infty}(\bftheta^*)$ with some but not necessarily all eigenvalues bounded; see \cite{CROWDER_1976}. In the sequel we, therefore, extend the definition suggested by \cite{PAZMAN_2007}.

Let $\{\bftau^{(n)}\}_{n\geq m}$ be a sequence of designs on $\calT_{n,\calD}$ such that $t_1^{(n)}=T_*$, $t_n^{(n)}=T^*$ and $\|\bftau^{(n)}\|\to0$, and let $\bftheta=(\bftheta_\romanI\T,\bftheta_\romanII\T)\T$ be a partition of the unknown parameter. Then the Fisher information matrix $\FIM_\romanI(\bftau^{(n)},\bftheta^*)$ corresponding to $\bftheta_\romanI$ is given by the Schur complement, for which we have
\begin{eqnarray*}
\FIM_\romanI(\bftau^{(n)},\bftheta^*) &=& \{\FIM(\bftau^{(n)},\bftheta^*)\}_{\bftheta_\romanI\bftheta_\romanI} -\{\FIM(\bftau^{(n)},\bftheta^*)\}_{\bftheta_\romanI\bftheta_\romanII}\{\FIM(\bftau^{(n)},\bftheta^*)\}_{\bftheta_\romanII\bftheta_\romanII}^{-1}\{\FIM(\bftau^{(n)},\bftheta^*)\}_{\bftheta_\romanII\bftheta_\romanI}\\
&\LoewnerDominatedBy& \{\FIM(\bftau^{(n)},\bftheta^*)\}_{\bftheta_\romanI\bftheta_\romanI} \nearrow \{\FIM_{\infty}(\bftheta^*)\}_{\bftheta_\romanI\bftheta_\romanI},
\end{eqnarray*}
as $n\to\infty$.
Additionally, the Loewner isotonicity of the information functions yields
\begin{equation*}
\Phi[\FIM_\romanI(\bftau^{(n)},\bftheta^*)]\leq\Phi[\{\FIM_{\infty}(\bftheta^*)\}_{\bftheta_\romanI\bftheta_\romanI}].
\end{equation*}

\begin{definition}\label{def:ULTIMATE_EFFICIENCY}
Let $\{\bftau^{(n)}\}_{n\geq m}$ be a sequence of designs on $\calT_{n,\calD}$ such that $t_1^{(n)}=T_*$, $t_n^{(n)}=T^*$ and $\|\bftau^{(n)}\|\to0$ as $n\to\infty$, and let $\bftheta=(\bftheta_\romanI\T,\bftheta_\romanII\T)\T$ be a partition of the unknown parameter, for which
\begin{equation*}
\lambda_{\max}\left(\{\FIM_{\infty}(\bftheta^*)\}_{\bftheta_\romanI\bftheta_\romanI}\right)<\infty
\end{equation*}
and
\begin{equation}\label{eqn:CONVERGENCE_CONDITION_FOR_ULTIMATE_EFFICIENCY}
\{\FIM(\bftau^{(n)},\bftheta^*)\}_{\bftheta_\romanI\bftheta_\romanII}\{\FIM(\bftau^{(n)},\bftheta^*)\}_{\bftheta_\romanII\bftheta_\romanII}^{-1}\{\FIM(\bftau^{(n)},\bftheta^*)\}_{\bftheta_\romanII\bftheta_\romanI} \to \0_{\dim(\bftheta_\romanI)\times\dim(\bftheta_\romanI)},
\end{equation}
as $n\to\infty$.
Then, the (local) ultimate efficiency of a design $\bftau\in\calT_{n,\calD}$ with respect to an information function $\Phi$, shortly (local) ultimate $\Phi$-efficiency, is the ratio
\begin{equation*}
\efficiency(\bftau\mid\Phi,\bftheta^*) = \frac{\Phi[\FIM_\romanI(\bftau,\bftheta^*)]}{\Phi[\{\FIM_{\infty}(\bftheta^*)\}_{\bftheta_\romanI\bftheta_\romanI}]}.
\end{equation*}
\end{definition}

Condition \eqref{eqn:CONVERGENCE_CONDITION_FOR_ULTIMATE_EFFICIENCY} in the definition of ultimate efficiency ensures that $\Phi[\FIM_\romanI(\bftau^{(n)},\bftheta^*)]\to\Phi[\{\FIM_{\infty}(\bftheta^*)\}_{\bftheta_\romanI\bftheta_\romanI}]$ as $n\to\infty$ and $\|\bftau^{(n)}\|\to0$.

\begin{theorem}\label{thm:ASYMPTOTIC_INFORMATION_MATRIX}
Let $\{X(t)\}_{t\geq0}$ be a process governed by stochastic differential equation \eqref{eqn:MODEL}, and let $\{\bftau^{(n)}\}_{n\geq m}$ be a sequence of designs on $\calT_{n,\calD}$ such that $t_1^{(n)}=T_*$, $t_n^{(n)}=T^*$ and $\|\bftau^{(n)}\|\to0$. Then
\begin{equation*}
\lim_{n\to\infty}\left(\FIM(\bftau^{(n)}_{},\bftheta^*)-\FIM_{\infty}(\bftheta^*)-\TheRest(\bftau^{(n)},\bftheta^*)\right)=\0_{m\times m},
\end{equation*}
where
\begin{eqnarray}
\FIM_{\infty}(\bftheta^*) &=& \frac{\frac{\partial \E[X(T_*)]}{\partial\bftheta}\frac{\partial \E[X(T_*)]}{\partial\bftheta\T}}{\Var[X(T_*)]}\Bigg|_{\bftheta^*} + \frac{1}{2}\frac{\partial \ln\Var[X(T_*)]}{\partial \bftheta}\frac{\partial \ln\Var[X(T_*)]}{\partial \bftheta\T}\Bigg|_{\bftheta^*}\label{eqn:ASYMPTOTIC_INFORMATION_MATRIX}\\
&&+\int_{T_*}^{T^*}\frac{\left.\frac{\partial f(t,x)}{\partial \bftheta}\right|_{x=\E[X(t)]}\left.\frac{\partial f(t,x)}{\partial \bftheta\T}\right|_{x=\E[X(t)]}+\frac{\partial b(t)}{\partial\bftheta}\frac{\partial b(t)}{\partial\bftheta\T}\Var[X(t)]}{\sigma^2(t)}\Bigg|_{\bftheta^*}\dif t\nonumber
\end{eqnarray}
and $\TheRest(\bftau^{(n)},\bftheta^*)=\frac{1}{2}\sum_{i=2}^n \frac{\partial\ln\sigma^2(t_i^{(n)})}{\partial\bftheta}\frac{\partial\ln\sigma^2(t_i^{(n)})}{\partial\bftheta\T}\big|_{\bftheta^*}$.
\end{theorem}

\begin{proof}
For the sake of simplicity we consider all functions and their derivatives to be evaluated at $\bftheta^*$. To prove the theorem we use Lemma~\ref{lem:INFORMATION_FOR_MARKOV_PROCESS}. Clearly,
\begin{equation*}
\FIM_{X(T_*)\mid X(0)}(\bftheta^*) = \frac{\frac{\partial \E[X(T_*)]}{\partial\bftheta}\frac{\partial \E[X(T_*)]}{\partial\bftheta\T}}{\Var[X(T_*)]} + \frac{1}{2}\frac{\partial \ln\Var[X(T_*)]}{\partial \bftheta}\frac{\partial \ln\Var[X(T_*)]}{\partial \bftheta\T}
\end{equation*}
gives the first two summands in \eqref{eqn:ASYMPTOTIC_INFORMATION_MATRIX}. Let $\bftau^{(n)}\in\calT_{n,\calD}$. We have
\begin{eqnarray}
\FIM_{X(t_i^{(n)})\mid X(t_{i-1}^{(n)})}(\bftheta^*) &=& \frac{1}{\Var[X(t_i^{(n)})\mid X(t_{i-1}^{(n)})]}\frac{\partial \E[X(t_i^{(n)})\mid X(t_{i-1}^{(n)})]}{\partial \bftheta}\frac{\partial \E[X(t_i^{(n)})\mid X(t_{i-1}^{(n)})]}{\partial \bftheta\T}\label{eqn:PARTIAL_FISHER_INFORMATION}\\
&&+\frac{1}{2}\frac{\partial \ln\Var[X(t_i^{(n)})\mid X(t_{i-1}^{(n)})]}{\partial \bftheta}\frac{\partial \ln\Var[X(t_i^{(n)})\mid X(t_{i-1}^{(n)})]}{\partial \bftheta\T},\nonumber
\end{eqnarray}
where
\begin{eqnarray*}
\E[X(t_i^{(n)})\mid X(t_{i-1}^{(n)})] &=& \Euler^{B(t_i^{(n)})-B(t_{i-1}^{(n)})}X(t_{i-1}^{(n)}) + \int_{t_{i-1}^{(n)}}^{t_i^{(n)}}\Euler^{B(t_i^{(n)})-B(\nu)}a(\nu)\dif\nu,\\
\Var[X(t_i^{(n)})\mid X(t_{i-1}^{(n)})] &=& \int_{t_{i-1}^{(n)}}^{t_i^{(n)}} \Euler^{2[B(t_i^{(n)})-B(\nu)]}\sigma^2(\nu)\dif\nu.
\end{eqnarray*}
From the Taylor expansion for $n\to\infty$ we obtain that
\begin{equation}\label{eqn:CONDITIONED_VARIANCE_EXPANSION}
\Var[X(t_i^{(n)})\mid X(t_{i-1}^{(n)})] = \sigma^2(t_i^{(n)})\Delta_i^{(n)} + \ord(\Delta_i^{(n)}),
\end{equation}
where $\Delta_i^{(n)}=t_i^{(n)}-t_{i-1}^{(n)}$ and $\ord(\Delta)/\Delta\to0$ as $\Delta\to0$. Consequently, for $n\to\infty$
\begin{equation*}
\frac{\partial \ln\Var[X(t_i^{(n)})\mid X(t_{i-1}^{(n)})]}{\partial \bftheta} \to \frac{\partial \ln\sigma^2(t_i^{(n)})}{\partial \bftheta}.
\end{equation*}
The only task left is to investigate the expectation of the first summand in \eqref{eqn:PARTIAL_FISHER_INFORMATION} with respect to $X(t_{i-1})$. Since $\E[X(t_i^{(n)})\mid X(t_{i-1}^{(n)})] = \E[X(t_i^{(n)})] + \Euler^{B(t_i^{(n)})-B(t_{i-1}^{(n)})}(X(t_{i-1}^{(n)}) - \E[X(t_{i-1}^{(n)})])$, we can write
\begin{multline*}
\E_{X(t_{i-1}^{(n)})}\left[\frac{\partial \E[X(t_i^{(n)})\mid X(t_{i-1}^{(n)})]}{\partial \bftheta}\frac{\partial \E[X(t_i^{(n)})\mid X(t_{i-1}^{(n)})]}{\partial \bftheta\T}\right] =\\ \frac{\partial \E[X(t_i^{(n)})]}{\partial \bftheta}\frac{\partial \E[X(t_i^{(n)})]}{\partial \bftheta\T} -\Euler^{B(t_i^{(n)})-B(t_{i-1}^{(n)})}\frac{\partial \E[X(t_i^{(n)})]}{\partial \bftheta}\frac{\partial \E[X(t_{i-1}^{(n)})]}{\partial \bftheta\T} -\Euler^{B(t_i^{(n)})-B(t_{i-1}^{(n)})}\frac{\partial \E[X(t_{i-1}^{(n)})]}{\partial \bftheta}\frac{\partial \E[X(t_i^{(n)})]}{\partial \bftheta\T}\\
+\Euler^{2[B(t_i^{(n)})-B(t_{i-1}^{(n)})]}\frac{\partial \E[X(t_{i-1}^{(n)})]}{\partial \bftheta}\frac{\partial \E[X(t_{i-1}^{(n)})]}{\partial \bftheta\T} +\frac{\partial [B(t_i^{(n)})-B(t_{i-1}^{(n)})]}{\partial \bftheta}\Var[X(t_{i-1}^{(n)})]\frac{\partial [B(t_i^{(n)})-B(t_{i-1}^{(n)})]}{\partial \bftheta\T}.
\end{multline*}
If we realize that
\begin{eqnarray*}
\Euler^{B(t_i^{(n)})-B(t_{i-1}^{(n)})}\frac{\partial \E[X(t_{i-1}^{(n)})]}{\partial \bftheta} &=& \frac{\partial \E[X(t_i^{(n)})]}{\partial \bftheta} - \frac{\partial}{\partial \bftheta}\int_{t_{i-1}^{(n)}}^{t_i^{(n)}}\Euler^{B(t_i^{(n)})-B(\nu)}a(\nu)\dif\nu -\frac{\partial [B(t_i^{(n)})-B(t_{i-1}^{(n)})]}{\partial \bftheta}\E[X(t_i^{(n)})]\\ && +\frac{\partial [B(t_i^{(n)})-B(t_{i-1}^{(n)})]}{\partial \bftheta}\int_{t_{i-1}^{(n)}}^{t_i^{(n)}}\Euler^{B(t_i^{(n)})-B(\nu)}a(\nu)\dif\nu,
\end{eqnarray*}
then some algebraic manipulation results in
\begin{multline*}
\E_{X(t_{i-1}^{(n)})}\left[\frac{\partial \E[X(t_i^{(n)})\mid X(t_{i-1}^{(n)})]}{\partial \bftheta}\frac{\partial \E[X(t_i^{(n)})\mid X(t_{i-1}^{(n)})]}{\partial \bftheta\T}\right] =\\ \left(\frac{\partial}{\partial \bftheta}\int_{t_{i-1}^{(n)}}^{t_i^{(n)}}\Euler^{B(t_i^{(n)})-B(\nu)}a(\nu)\dif\nu+\frac{\partial [B(t_i^{(n)})-B(t_{i-1}^{(n)})]}{\partial \bftheta}\E[X(t_i^{(n)})] -\frac{\partial [B(t_i^{(n)})-B(t_{i-1}^{(n)})]}{\partial \bftheta}\int_{t_{i-1}^{(n)}}^{t_i^{(n)}}\Euler^{B(t_i^{(n)})-B(\nu)}a(\nu)\dif\nu\right)\\
\times\left(\frac{\partial}{\partial \bftheta\T}\int_{t_{i-1}^{(n)}}^{t_i^{(n)}}\Euler^{B(t_i^{(n)})-B(\nu)}a(\nu)\dif\nu+\frac{\partial [B(t_i^{(n)})-B(t_{i-1}^{(n)})]}{\partial \bftheta\T}\E[X(t_i^{(n)})] -\frac{\partial [B(t_i^{(n)})-B(t_{i-1}^{(n)})]}{\partial \bftheta\T}\int_{t_{i-1}^{(n)}}^{t_i^{(n)}}\Euler^{B(t_i^{(n)})-B(\nu)}a(\nu)\dif\nu\right)\\
+\frac{\partial [B(t_i^{(n)})-B(t_{i-1}^{(n)})]}{\partial \bftheta}\frac{\partial [B(t_i^{(n)})-B(t_{i-1}^{(n)})]}{\partial \bftheta\T}\Var[X(t_{i-1}^{(n)})].
\end{multline*}
From the fact $\Var[X(t_{i-1}^{(n)})] = \Var[X(t_i^{(n)})] - [2b(t_i^{(n)})\Var[X(t_i^{(n)})]+\sigma^2(t_{i-1}^{(n)})]\Delta_i^{(n)}$ we get the Taylor expansion
\begin{multline*}
\E_{X(t_{i-1}^{(n)})}\left[\frac{\partial \E[X(t_i^{(n)})\mid X(t_{i-1}^{(n)})]}{\partial \bftheta}\frac{\partial \E[X(t_i^{(n)})\mid X(t_{i-1}^{(n)})]}{\partial \bftheta\T}\right] =\\
\left(\frac{\partial a(t_i^{(n)})}{\partial \bftheta} + \frac{\partial b(t_i^{(n)})}{\partial \bftheta}\E[X(t_i^{(n)})]\right)\left(\frac{\partial a(t_i^{(n)})}{\partial \bftheta\T} + \frac{\partial b(t_i^{(n)})}{\partial \bftheta\T}\E[X(t_i^{(n)})]\right)(\Delta_i^{(n)})^2\\
+ \frac{\partial b(t_i^{(n)})}{\partial \bftheta}\frac{\partial b(t_i^{(n)})}{\partial \bftheta\T}\Var[X(t_i^{(n)})](\Delta_i^{(n)})^2 + \ord\left((\Delta_i^{(n)})_{}^2\right),
\end{multline*}
which together with expectation of the conditioned Fisher information matrix \eqref{eqn:PARTIAL_FISHER_INFORMATION} with respect to $X(t_{i-1})$ and formula \eqref{eqn:CONDITIONED_VARIANCE_EXPANSION} for variance yields the statement of the theorem.
\end{proof}

Because the underlying model considers $\theta_3$ to be the only parameter of the volatility function $\sigma(t)$, all elements of 
$\TheRest(\bftau^{(n)},\bftheta^*)$ are zero except $\{\TheRest(\bftau^{(n)},\bftheta^*)\}_{\theta_3\theta_3}$, which tends to infinity. All elements of the asymptotic Fisher information matrix $\FIM_{\infty}(\bftheta^*)+\TheRest(\bftau^{(n)},\bftheta^*)$ are, therefore, bounded for $n\to\infty$ as well, except the diagonal entry $\{\FIM_{\infty}(\bftheta^*)+\TheRest(\bftau^{(n)},\bftheta^*)\}_{\theta_3\theta_3}$. That is, if we denote $\bftheta_\romanI = (X_0,\bftheta_1\T,\bftheta_2\T)\T$ then
\begin{equation*}
\{\FIM_{\infty}(\bftheta^*)+\TheRest(\bftau^{(n)},\bftheta^*)\}_{\bftheta_\romanI\theta_3}\{\FIM_{\infty}(\bftheta^*)+\TheRest(\bftau^{(n)},\bftheta^*)\}_{\theta_3\theta_3}^{-1}\{\FIM_{\infty}(\bftheta^*)+\TheRest(\bftau^{(n)},\bftheta^*)\}_{\theta_3\bftheta_\romanI}\to0,\;\text{as}\;n\to\infty,
\end{equation*}
so for a parameter subvector consisting of $X_0$, $\bftheta_1$ and $\bftheta_2$ we can, in line with Definition \ref{def:ULTIMATE_EFFICIENCY}, compute the ultimate $\Phi$-efficiency of designs.

Besides the form of the asymptotic information matrix, Theorem~\ref{thm:ASYMPTOTIC_INFORMATION_MATRIX} provides a basis for analysis of the maximum likelihood estimator consistency. The structure of the Fisher information matrix $\FIM_{\infty}(\bftheta^*)+\TheRest(\bftau^{(n)},\bftheta^*)$ indicates that the only unknown parameter we can estimate consistently is $\theta_3$. This phenomenon has a natural explanation: Unlike $X_0$, $\bftheta_1$ or $\bftheta_2$, in equation \eqref{eqn:MODEL} the parameter $\theta_3$ is connected with the differential of the Wiener process $\{W(t)\}_{t\geq0}$, which is characterized by the fractal property called Brownian scaling. Consequently, the more observations we perform the more information about $\theta_3$ we gain. An analogous formulation of this result can be found in the literature on the stochastic differential equations, see, e.g., \cite{IACUS_2008} for a brief survey, but is also noted in selected papers on inference in regression problems with correlated observations; see \cite{PAZMAN_1965}.

\section{Processes of Ornstein-Uhlenbeck type}\label{sec:PROCESSES_OF_ORNSTEIN_UHLENBECK_TYPE}

In the previous sections we presented an analysis of a generalised form of Ornstein-Uhlenbeck process described by stochastic differential equation \eqref{eqn:MODEL}, which enables us to evaluate a quality of sampling designs. A reasonable question is whether we can use the results also for other processes than those governed by equation \eqref{eqn:MODEL}.

The motivation comes from the Fisher-Neymann factorization theorem; see, for instance, Theorem~6.5 in \cite{LEHMANN_CASELLA_1998}. More precisely, if we apply a sufficient statistic to measurements then the Fisher information matrix remains unchanged. Henceforth, we can define the following class of stochastic differential equations.

\begin{definition}
Let $\mu(t,y)$ and $\gamma(t,y)$ be sufficiently smooth, and let a process $\{Y(t)\}_{t\geq0}$ be governed by a stochastic differential equation
\begin{equation*}
\dif Y(t) = \mu(t,Y(t))\dif t + \gamma(t,Y(t))\dif W(t).
\end{equation*}
If there exist sufficiently smooth functions $a(t)$, $b(t)$, $\sigma(t)$ and $\varphi(t,y)$, where $\varphi$ is bijective in $y$ and $\frac{\partial \varphi(t,y)}{\partial \bftheta}=0$ for all $t\geq0$ and $y\in\Reals$, such that a process $\{X(t)\}_{t\geq0}=\{\varphi(t,Y(t))\}_{t\geq0}$ is governed by an equation
\begin{equation*}
\dif X(t) = [a(t)+b(t)X(t)]\dif t + \sigma(t)\dif W(t),
\end{equation*}
then we say that the process $\{Y(t)\}_{\geq0}$ is a process of Ornstein-Uhlenbeck type with associated coefficients $a(t)$, $b(t)$ and $\sigma(t)$. We denote this fact by $\{Y(t)\}\in\OUPClass_{a(t),b(t),\sigma(t)}$.
\end{definition}

A candidate for the desired sufficient statistic, which transforms the volatility of the original process to the volatility of a desired form, is
\begin{equation*}
\varphi(t,y) = \int \frac{\sigma(t)}{\gamma(t,y)}\dif y.
\end{equation*}
The condition that $\frac{\partial \varphi(t,y)}{\partial \bftheta}=0$ for all $t\geq0$ and $y\in\Reals$ might not be easy to verify in advance, because we do not know the form of $\sigma(t)$. Nevertheless, if we can write the diffusion term in a separable form $\gamma(t,y)=\sigma(t)g(t,y)$, where $\frac{\partial g(t,y)}{\partial \bftheta}=0$ for all $t\geq0$ and $y\in\Reals$, then $\varphi$ is a sufficient statistic. We remark that $\varphi$ depends on a reciprocal of $\gamma(t,y)$, thus we might need to impose further positivity conditions on the domain interior of $Y(t)$ for all $t$, which is out of scope of the presented paper.

In the sequel we propose a way how to verify whether a given process is of Ornstein-Uhlenbeck type. Let $\psi(t,y)$ be the inverse function to $\varphi(t,y)$, that is, $\psi(t,\varphi(t,y))=y$. Then $\{Y(t)\}_{t\geq0}=\{\psi(t,X(t))\}_{t\geq0}$. It\=o's lemma implies
\begin{equation}\label{eqn:PDE_FOR_COEFFICIENTS}
\left.\frac{\partial \psi}{\partial t}\right|_{x=\varphi(t,y)} + \left.\frac{\partial \psi}{\partial x}\right|_{x=\varphi(t,y)}f(t,\varphi(t,y)) + \frac{1}{2}\left.\frac{\partial^2 \psi}{\partial x^2}\right|_{x=\varphi(t,y)}\sigma^2(t) = \mu(t,y).
\end{equation}
By substituting the relations for inverse functions
\begin{eqnarray*}
\left.\frac{\partial \psi}{\partial t}\right|_{x=\varphi(t,y)}= -\frac{\partial \varphi}{\partial t}/\frac{\partial \varphi}{\partial y},\;\left.\frac{\partial \psi}{\partial x}\right|_{x=\varphi(t,y)}=-1/\frac{\partial \varphi}{\partial y},\;\left.\frac{\partial^2 \psi}{\partial x^2}\right|_{x=\varphi(t,y)}=-\frac{\partial^2 \varphi}{\partial y^2}/\left(\frac{\partial \varphi}{\partial y}\right)^3.\\
\end{eqnarray*}
to equation \eqref{eqn:PDE_FOR_COEFFICIENTS} we obtain the following theorem.

\begin{theorem}\label{thm:PROCESS_OF_ORNSTEIN_UHLENBECK_TYPE_CONDITIONS}
Let the process $\{Y(t)\}_{t\geq0}$ be driven by a stochastic differential equation
\begin{equation*}
\dif Y(t) = \mu(t,Y(t))\dif t + \sigma(t)g(t,Y(t))\dif W(t),
\end{equation*}
where the functions $\mu(t,y)$ and $g(t,y)$ are sufficiently smooth and $\frac{\partial g(t,y)}{\partial\bftheta}=0$ for all $t\geq0$ and $y\in\Reals$. If  the condition
\begin{equation}\label{eqn:PROCESS_OF_ORNSTEIN_UHLENBECK_TYPE_CONDITIONS}
\frac{\dif}{\dif y}\int\frac{\dif t}{g(t,y)} = a(t) +b(t)\int\frac{\dif y}{g(t,y)} +\frac{1}{2}\sigma^2(t)\frac{\partial g(t,y)}{\partial y}-\frac{\mu(t,y)}{g(t,y)}
\end{equation}
is satisfied for some functions $a(t)$ and $b(t)$ then $\{Y(t)\}_{t\geq0}\in\OUPClass_{a(t),b(t),\sigma(t)}$.
\end{theorem}

In some instances we consider a process governed by an autonomous stochastic differential equation of the form
\begin{equation}\label{eqn:AUTONOMOUS_MODEL}
\dif Y(t) = \mu(Y(t))\dif t + \sigma g(Y(t))\dif W(t),
\end{equation}
where the structure of the drift function $\mu(y)$ is usually based on an essential theory in the given research field. On the contrary, the choice of $g(y)$ might be artificial; we choose a diffusion that fits some arrangements. By a differentiation of equation \eqref{eqn:PROCESS_OF_ORNSTEIN_UHLENBECK_TYPE_CONDITIONS} with respect to $y$ we get

\begin{corollary}
Let the process $\{Y(t)\}_{t\geq0}\in\OUPClass_{a,b,\sigma}$ be driven by a stochastic differential equation \eqref{eqn:AUTONOMOUS_MODEL} and let $\mu(y)$ be given. Then $g(y)$ solves the ordinary differential equation
\begin{equation*}
\mu(y) + \left(b-\frac{\partial \mu(y)}{\partial y}\right)g(y) + \frac{1}{2}\sigma^2 g^2(y)\frac{\partial^2 g(y)}{\partial y^2} = 0.
\end{equation*}
\end{corollary}

If a positive solution $g(y)$ at least approximately corresponds to the experimental setting, we can use the proposed methodology for assessment of a design quality.

\section{Example: Gompertz model of tumour growth}\label{sec:GOMPERTZ_MODEL_OF_TUMOUR_GROWTH}

Tumour growth models play an important role in therapeutic guidance. \cite{GOMPERTZ_1825} proposed in his pioneering paper a growth model, which became a base for many studies in cancer research and various modifications of this growth law were introduced; we refer the reader to \cite{NORTON_1988} and \cite{SPEER_PETROSKY_RETSKY_ET_AL_1984} for a brief survey. Although the work of \cite{CAMERON_1997} demonstrates that the Gompertz model is not the best choice for the breast cancer, it seems to fit well data for multiple myeloma; see \cite{SULLIVAN_SALMON_1972}.

In the Gompertz growth model the expected size of the tumour $Y_{\rho,\delta}(t)$ solves an ordinary differential equation
\begin{equation}\label{eqn:DETERMINISTIC_GOMPERTZ_MODEL}
\frac{\dif Y(t)}{\dif t} = \rho Y(t) - \delta Y(t)\ln Y(t),\;Y(0)=Y_0\;\text{known},
\end{equation}
where the intrinsic growth rate $\rho$ and the growth deceleration factor $\delta$ are unknown parameters, and the initial size $Y_0$ is usually obtained from the first-time detection. The solution to equation \eqref{eqn:DETERMINISTIC_GOMPERTZ_MODEL} is known in an explicit closed form
\begin{equation}\label{eqn:GOMPERTZ_MODEL_RESPONSE_FUNCTION}
Y_{\rho,\delta}(t) = \Euler^{\Euler^{-\delta t}\ln Y_0 + (1-\Euler^{-\delta t})\frac{\rho}{\delta}},
\end{equation}
which might be a factor explaining its popularity in applications.

To improve the quality of statistical inference it is important to choose a proper experimental design. Locally D-optimal designs for model \eqref{eqn:DETERMINISTIC_GOMPERTZ_MODEL} were studied by \cite{LI_2012}, who assumed that the observations $\tilde{Y}_{ij}$ satisfy the regression model
\begin{equation*}
\tilde{Y}_{ij} = Y_{\rho,\delta}(t_i) + \varepsilon_{ij},\;i=1,\ldots,n,\;j=1,\ldots,k_i.
\end{equation*}
Here, $\tilde{Y}_{ij}$ are observations of tumour size, $t_i$'s are the design times, at which we perform $k_i$ replications, and $\varepsilon_{ij}$'s are independent homoscedastic errors.

Notice that the response function \eqref{eqn:GOMPERTZ_MODEL_RESPONSE_FUNCTION} is nonlinear and a construction of optimal designs using local linearization might require prior estimates of the true values of parameters $\rho$ and $\delta$. We can obtain prior estimates using, for instance, multiresponse regression models, which capture relations between estimates and different physiologic factors using the data from former experiments.

In this section we consider a variant of the Gompertz growth model, where the size of the tumour $\{Y(t)\}_{t\geq0}$ is described by a stochastic differential equation
\begin{equation}\label{eqn:GOMPERTZ_MODEL}
\dif Y(t) = [\rho Y(t) - \delta Y(t)\ln Y(t)]\dif t + \gamma Y(t)\dif W(t),\;Y(0)=Y_0\;\text{known},
\end{equation}
see \cite{LO_2009} for more details. Compared to the model of \cite{LI_2012}, for the process $\{Y(t)\}_{t\geq0}$ governed by equation \eqref{eqn:GOMPERTZ_MODEL} we cannot replicate the observations (unless we have more patients), and so the amount of information about parameters might be limited.

Remark that the structure of volatility in equation \eqref{eqn:GOMPERTZ_MODEL} implies small fluctuations in the tumour size for small tumours and, on the contrary, if the size of the tumour is large then the fluctuations in the size of the tumour are large, too.

If we define $\{X(t)\}_{t\geq0}=\{\ln Y(t)\}_{t\geq0}$ then $\{X(t)\}_{t\geq0}$ solves
\begin{equation*}
\dif X(t) = \left(\rho-\frac{1}{2}\gamma^2 - \delta X(t)\right)\dif t + \gamma\dif W(t).
\end{equation*}
That is, the process $\{Y(t)\}_{t\geq0}$ is a process of Ornstein-Uhlenbeck type with the associated coefficients
\begin{equation*}
a_{\rho,\gamma}(t)=\rho-\frac{1}{2}\gamma^2,\; b_{\delta}(t) = -\delta\; \text{and}\;\sigma_\gamma(t)=\gamma.
\end{equation*}
The parameter $\rho$ is of type $\bftheta_1$, $\delta$ is of the type $\bftheta_2$ and $\gamma$ is of the type $\theta_3$, cf. equation \eqref{eqn:MODEL}. By setting $\bftheta_\romanI=(\rho,\delta)\T$ get the asymptotic Fisher information matrix for $\bftheta_\romanI$
\begin{equation*}
\{\FIM_{\infty}(\bftheta^*)\}_{\bftheta_\romanI\bftheta_\romanI} = \frac{\frac{\partial \E[X(T_*)]}{\partial \bftheta_\romanI}\frac{\partial \E[X(T_*)]}{\partial \bftheta_\romanI\T}}{\Var[X(T_*)]} + \frac{\partial \ln\Var[X(T_*)]}{\partial \bftheta_\romanI}\frac{\partial \ln\Var[X(T_*)]}{\partial \bftheta_\romanI\T} +\frac{1}{\gamma^2}\int_{T_*}^{T^*}\begin{pmatrix}
1& -\E[X(t)]\\
-\E[X(t)]&\E^2[X(t)] + \Var[X(t)]
\end{pmatrix}\dif t,
\end{equation*}
where $\E[X(t)] = \Euler^{-\delta t}\ln Y_0 + \frac{1}{\delta}(\rho-\frac{1}{2}\gamma^2)(1-\Euler^{-\delta t})$, $\Var[X(t)] = \frac{\gamma^2}{2\delta}(1-\Euler^{-2\delta t})$, and the integration is meant to be component-wise. We note that with exception of the information contained in the first observation (which is zero if $T^*\to0$), the asymptotic information depends only on the mean value and the variance of the process.

\begin{figure}
\centering

\begin{minipage}[t]{.35\textwidth}
\centering
\includegraphics[width=\textwidth]{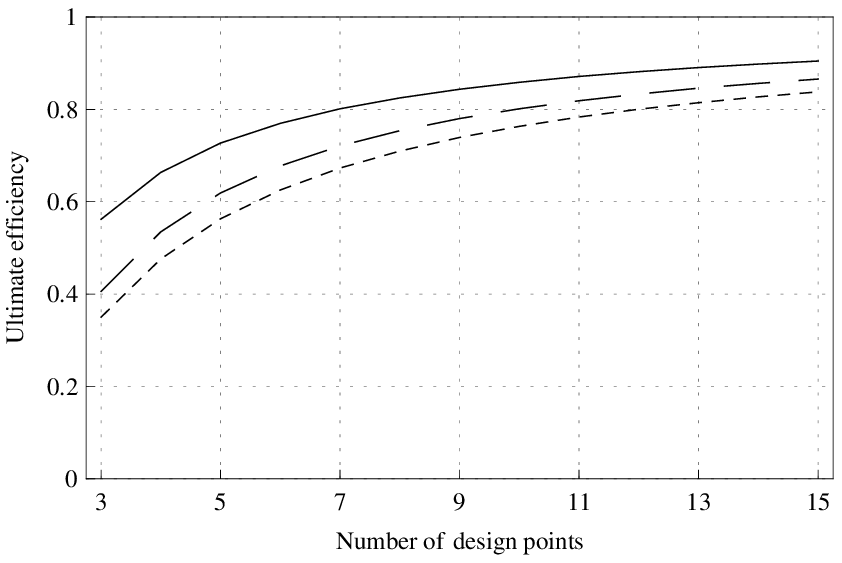}\\
\tiny a) $\rho^*=1$, $\delta^*=1$, $\gamma^*=1$
\end{minipage}
\begin{minipage}[t]{.35\textwidth}
\centering
\includegraphics[width=\textwidth]{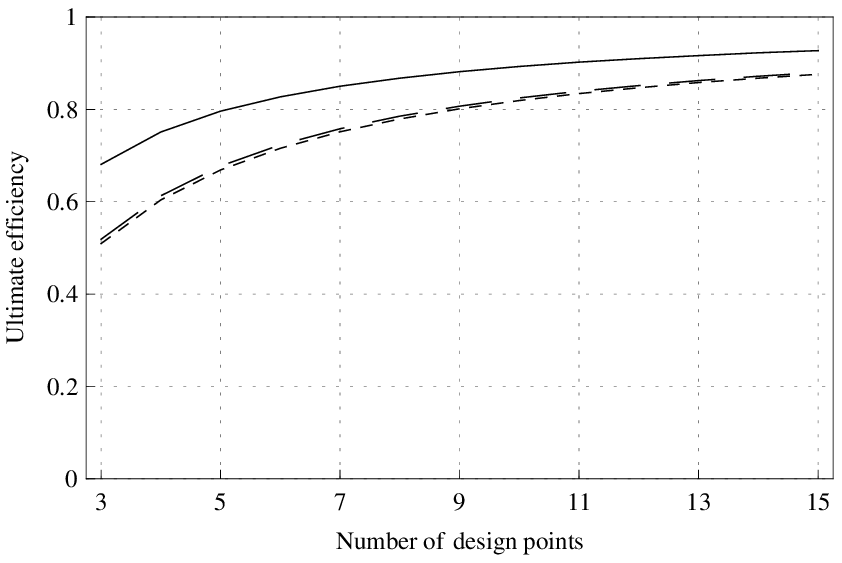}\\
\tiny b) $\rho^*=3$ ,$\delta^*=1$, $\gamma^*=1$
\end{minipage}

\bigskip

\begin{minipage}[t]{.35\textwidth}
\centering
\includegraphics[width=\textwidth]{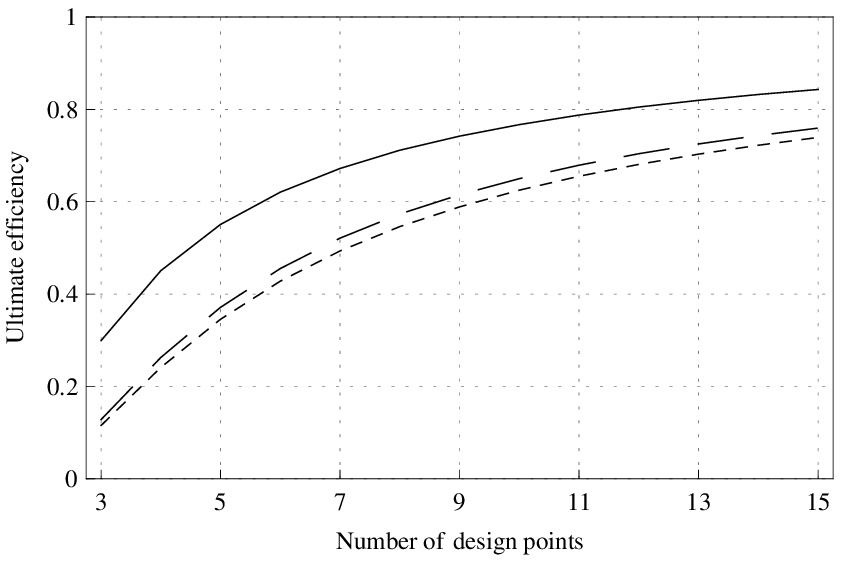}\\
\tiny c) $\rho^*=1$ ,$\delta^*=3$, $\gamma^*=1$
\end{minipage}
\begin{minipage}[t]{.35\textwidth}
\centering
\includegraphics[width=\textwidth]{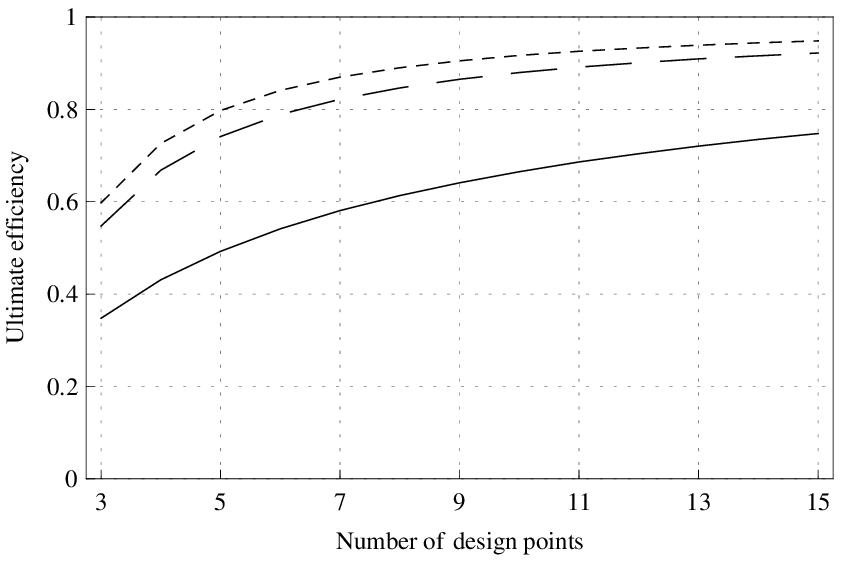}\\
\tiny d) $\rho^*=1$, $\delta^*=1$, $\gamma^*=3$
\end{minipage}

\caption{Ultimate efficiencies of equidistant $n$-point designs on experimental domain $\calD=[1,2]$ with $t_1=1$, $t_n=2$, $Y_0=1$ and different prior estimates $\rho^*$, $\delta^*$ and $\gamma^*$. We consider D-optimality criterion (------), E-optimality criterion (- - -) and A-optimality criterion (--- ---).}%
\label{fig:EFFICIENCIES}%
\end{figure}

Figure \ref{fig:EFFICIENCIES} depicts a dependence of the ultimate efficiency of an equidistant $n$-point design on the number of design points. In this example we consider $\calD=[1,2]$ and designs with $t_1=1$ and $t_n=2$, and different prior values of $\rho$, $\delta$ and $\gamma$. The nonlinearity of the regression model results in a contradiction with a generally accepted opinion that under correlated observations it is sufficient to perform a few trials to attain a very high level of ultimate efficiency. For example, to attain $80\%$ local ultimate efficiency for E-optimality if $\rho=1$, $\delta=1$ and $\gamma=3$ it is enough to perform $5$ trials (see Figure \ref{fig:EFFICIENCIES}d), and, on the other hand, if $\rho=1$, $\delta=3$, $\gamma=1$ then for $15$ observations the ultimate efficiency for E-optimality is slightly above $70\%$ (see Figure \ref{fig:EFFICIENCIES}c). From a practical point of view this fact says that if we have more objects to be studied in the experiment then we should not perform an  equal number of observations for each subject, but realocate the number of observations due to physiologic traits so that the overall experimental framework is efficient.

\subsection*{Acknowledgement}
The research was supported by the Slovak VEGA Grants No. 2/0038/12 and 1/0163/13, and by the Comenius University Grant No. UK/87/2012. The author would like to thank Prof. Andrej P\'azman and Assoc. Prof. Radoslav Harman for their comments that influenced the course of this paper, and to Prof. Daniel \v{S}ev\v{c}ovi\v{c} for consultations on specific topics.

\bibliographystyle{elsarticle-harv}
\bibliography{References}

\end{document}